\title{Paley and the Paley graphs}

\author{
Gareth~A.~Jones\\
}

\documentclass[12pt]{article}
\usepackage{color}
\usepackage[latin1]{inputenc}
\usepackage{a4wide}
\usepackage{latexsym}
\usepackage{amsmath}
\usepackage{amssymb}
\usepackage{amsfonts}
\usepackage{tikz}
\newtheorem{thm}{Theorem}[section]
\newtheorem{lemma}[thm]{Lemma}

\newcommand{\F}{\mathbb F}

\date{}
\begin{document}

\maketitle

\begin{abstract}
This paper discusses some aspects of the history of the Paley graphs and their automorphism groups.
\end{abstract}


\medskip

\section{Introduction}

Anyone who seriously studies algebraic graph theory or finite permutation groups will, sooner or later, come across the Paley graphs and their automorphism groups. The most frequently cited sources for these are respectively Paley's 1933 paper~\cite{Pal}, and Carlitz's 1960 paper~\cite{Car}. It is remarkable that neither of those papers uses the concepts of graphs, groups or automorphisms. Indeed, one cannot find these three terms, or any synonyms for them, in those papers: Paley's  paper is entirely about the construction of what are now called Hadamard matrices, while Carlitz's is entirely about permutations of finite fields.

The aim of the present paper is to explore how this strange situation came about, by explaining the background to these two papers and how they became associated with the Paley graphs. This involves describing various links with other branches of mathematics, such as matrix theory, number theory, design theory, coding theory, finite geometry, polytope theory and group theory. The paper is organised in two main parts, the first covering the graphs and the second their automorphism groups, each largely in historical order. However, in order to establish basic concepts we start with the definition and elementary properties of the Paley graphs.


\section{Definition and properties of the Paley graphs}

The Paley graph $P(q)$ has vertex set $V=\F:=\F_q$, a field of prime power order $q=p^e\equiv 1$ mod~$(4)$, with two vertices $u$ and $v$ adjacent if and only if $u-v$ is an element of the set
\[S=\{x^2\mid x\in\F,\, x\ne 0\}\]
of quadratic residues (non-zero squares) in $\F$. It is thus the Cayley graph~\cite{Cay} for the additive group of $\F$, with $S$ as the connection set.  The choice of $q$ ensures that $-1\in S$, so $S=-S$ and $P(q)$ is an undirected graph; the fact that $S$ generates the additive  group ensures that $P(q)$ is connected. The neighbours of a vertex $v$ are the elements of $S+v$, so its valency is $|S|=(q-1)/2$. (If $q\equiv 3$ mod~$(4)$ this construction gives a directed graph, in fact a tournament, since each pair of vertices $u\ne v$ are joined by a unique arc $u\to v$, where $v-u\in S$.)

\begin{figure}[h!]
\begin{center}
\begin{tikzpicture}[scale=0.2, inner sep=0.8mm]

\draw [dashed] (15,15) to (-15,15) to (-15,-15) to (15,-15) to (15,15);

\draw [thick] (-15,10) to (15,10);
\draw [thick] (-15,0) to (15,0);
\draw [thick] (-15,-10) to (15,-10);

\draw [thick] (-10,-15) to (-10,15);
\draw [thick] (0,-15) to (0,15);
\draw [thick] (10,-15) to (10,15);

\node (a) at (-10,10) [shape=circle, fill=black] {};
\node (b) at (0,10) [shape=circle, fill=black] {};
\node (c) at (10,10) [shape=circle, fill=black] {};
\node (d) at (-10,0) [shape=circle, fill=black] {};
\node (e) at (0,0) [shape=circle, fill=black] {};
\node (f) at (10,0) [shape=circle, fill=black] {};
\node (g) at (-10,-10) [shape=circle, fill=black] {};
\node (h) at (0,-10) [shape=circle, fill=black] {};
\node (i) at (10,-10) [shape=circle, fill=black] {};

\node at (-6.5,8) {$-1+i$};
\node at (2,8) {$i$};
\node at (12.5,8) {$1+i$};
\node at (-8,-2) {$-1$};
\node at (2,-2) {$0$};
\node at (12,-2) {$1$};
\node at (-6.5,-12) {$-1-i$};
\node at (2,-12) {$-i$};
\node at (12.5,-12) {$1-i$};

\end{tikzpicture}

\end{center}
\caption{$P(9)$ drawn on a torus}
\label{P(9)}
\end{figure}
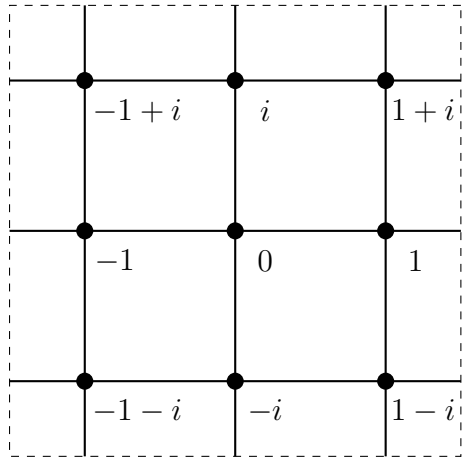

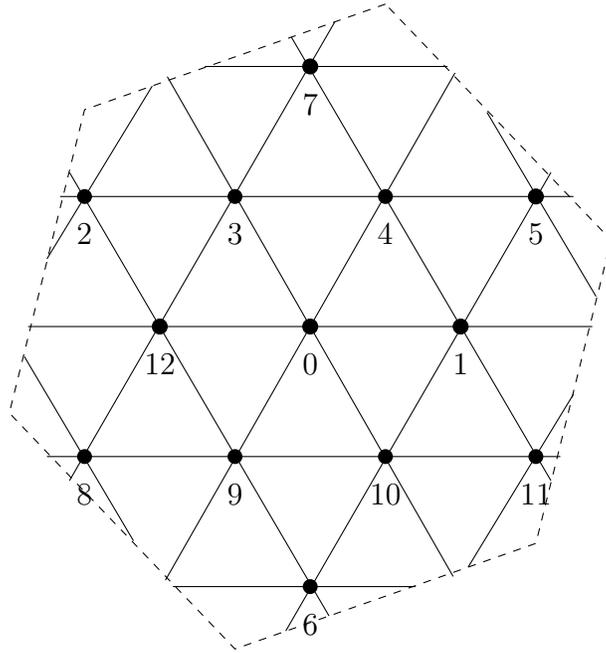
\begin{figure}[h!]
\label{cube}
\begin{center}
\begin{tikzpicture}[scale=0.5, inner sep=0.7mm]

\node (a) at (0,0) [shape=circle, draw, fill=black] {}; 
\node (b) at (4,0) [shape=circle, draw, fill=black] {}; 
\node (d) at (2,3.46) [shape=circle, fill=black] {};  
\node (c) at (-2,3.46) [shape=circle, fill=black] {};  
\node (g) at (-4,0) [shape=circle, draw, fill=black] {}; 
\node (e) at (-2,-3.46) [shape=circle, fill=black] {};  
\node (f) at (2,-3.46) [shape=circle, fill=black] {};  

\node at (0,-1) {$0$};
\node at (4,-1) {$1$};
\node at (2,2.46) {$4$};
\node at (-2,2.46) {$3$};
\node at (-4,-1) {$12$};
\node at (-2,-4.46) {$9$};
\node at (2,-4.46) {$10$};

\node (u) at (6,3.46) [shape=circle, draw, fill=black] {}; 
\node (v) at (0,6.92) [shape=circle, draw, fill=black] {}; 
\node (w) at (-6,3.46) [shape=circle, fill=black] {};  
\node (x) at (-6,-3.46) [shape=circle, fill=black] {};  
\node (y) at (0,-6.92) [shape=circle, fill=black] {};  
\node (z) at (6,-3.46) [shape=circle, fill=black] {};  

\node at (6,2.46) {$5$};
\node at (0,5.92) {$7$};
\node at (-6,2.46) {$2$};
\node at (-6,-4.46) {$8$};
\node at (0,-7.92) {$6$};
\node at (6,-4.46) {$11$};

\draw [dashed] (2,8.58) to (-6,5.77 ) to (-8 , -2.3) to (-2,-8.58) to (6,-5.77 ) to (8,2.3) to (2,8.58);

\draw (7.5,0) to (-7.5,0); %
\draw ( 7,3.46) to (-6.65, 3.46); %
\draw ( -7,-3.46) to (6.65,-3.46); %
\draw ( -2.8,6.92) to (3.65,6.92); %
\draw ( 2.8,-6.92) to (-3.65,-6.92); %
\draw (3.85,6.75) to (-3.85,-6.75); %
\draw (3.8,-6.65) to (-3.8,6.65); %

\draw (-0.5,7.7) to (6.4,-4.1); %
\draw (0.5,-7.7) to (-6.4,4.1);  %
\draw (0.65,8.1) to (-6.4,-4.1);  %
\draw (-0.65,-8.1) to (6.4,4.1); %
\draw (4.7,5.7) to (7.6,0.8); %
\draw (-4.7,-5.7) to (-7.6,-0.8); %
\draw (-4.2,6.4) to (-7,1.8); 
\draw (4.2,-6.4) to (7,-1.8); 
 
  \end{tikzpicture}

\end{center}
\caption{$P(13)$ drawn on a torus }
\label{P(13)torus}
\end{figure}

For example, Figure~\ref{P(9)} shows $P(9)$ drawn on a torus, formed by identifying opposite sides of the outer square. Here $\F_9=\F_3[i]$ where $i^2=-1$; in other words, this map is the quotient of a Cayley map for the additive group ${\mathbb Z}[i]$ of Gaussian integers, modulo the ideal $(3)$. In the notation of Coxeter and Moser~\cite[Ch.~8]{CM}, this is the map $\{4,4\}_{3,0}$. There is an analogous chiral pair of torus embeddings of $P(13)$ as the triangular maps $\{3,6\}_{3,1}$ and $\{3,6\}_{1,3}$; Figure~\ref{P(13)torus} shows the former, with opposite sides of the outer hexagon identified.

\begin{figure}[h!]
\begin{center}
\begin{tikzpicture}[scale=0.25, inner sep=0.8mm]

\node (0) at (0,10) [shape=circle, fill=black] {};
\node (1) at (4.65,8.85) [shape=circle, fill=black] {};
\node (2) at (8.24,5.67) [shape=circle, fill=black] {};
\node (3) at (9.93,1.20) [shape=circle, fill=black] {};
\node (4) at (9.35,-3.53) [shape=circle, fill=black] {};
\node (5) at (6.53,-7.48) [shape=circle, fill=black] {};
\node (6) at (2.39,-9.71) [shape=circle, fill=black] {};
\node (7) at (-2.39,-9.71) [shape=circle, fill=black] {};
\node (8) at (-6.53,-7.48) [shape=circle, fill=black] {};
\node (9) at (-9.35,-3.53) [shape=circle, fill=black] {};
\node (10) at (-9.93,1.20) [shape=circle, fill=black] {};
\node (11) at (-8.24,5.67) [shape=circle, fill=black] {};
\node (12) at (-4.65,8.85) [shape=circle, fill=black] {};

\draw [thick] (0) to (1) to (2) to (3) to (4) to (5) to (6) to (7) to (8) to (9) to (10) to (11) to (12) to (0);
\draw [thick, blue] (0) to (3) to (6) to (9) to (12) to (2) to (5) to (8) to (11) to (1) to (4) to (7) to (10) to (0);
\draw [thick, red] (0) to (4) to (8) to (12) to (3) to (7) to (11) to (2) to (6) to (10) to (1) to (5) to (9) to (0);

\end{tikzpicture}

\end{center}
\caption{$P(13)$, showing dihedral symmetry}
\label{P(13)}
\end{figure}
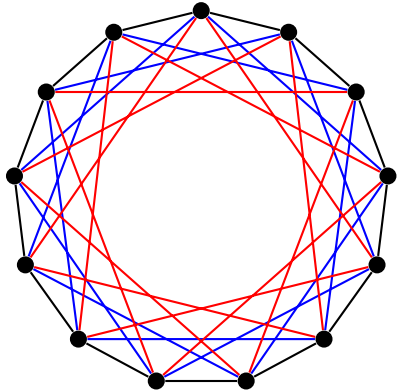

Figure~\ref{P(13)} shows $P(13)$, exhibiting its dihedral symmetry under the automorphisms $v\mapsto \pm v+b$, $b\in\F_{13}$. Vertices are identified with $0, 1, \ldots, 12$ in cyclic order, and edges $uv$ are coloured black, blue or red as $u-v=\pm 1, \pm 3$ or $\pm 4$ respectively.

It is clear from the definition that further combinatorial properties of the graphs $P(q)$ will depend on the properties of quadratic residues in finite fields. For any odd prime power $q$, let $Q$ be the {\sl Jacobsthal matrix\/} for $\F=\F_q$: this has rows and columns indexed by the elements of $\F$, with $(u,v)$ entry $\chi(v-u)$ for all $u, v\in\F$, where $\chi:\F\to{\mathbb C}$ is the {\sl quadratic residue character\/} of $\F$, defined by
\[\chi(x)=
\begin{cases}
0 & {\rm if} \; x=0,\\
1 & {\rm if} \; x\in S,\\
-1 & {\rm otherwise}.
\end{cases}.\]
Thus the restriction of $\chi$ to the multiplicative group $\F^*:=\F\setminus\{0\}$ is a group epimorphism $\F^*\to\{\pm 1\}$ with kernel $S$. When $q$ is a prime $p$, $\chi(x)$ is the Legendre symbol $(\frac{x}{p})$. 

The matrix $Q$ is symmetric or skew-symmetric as $q\equiv 1$ or $3$ mod~$(4)$, with
\[QJ=JQ=0\quad{\rm and}\quad QQ^T=qI-J,\]
where $J$ is the matrix of order $q$ with all entries equal to $1$. The first two equations here are obvious. For the last equation, note that $QQ^T$ has $(u,v)$ entry
$\sum_{w\in\F}\chi(w-u)\chi(w-v)$ given by the following lemma:

\begin{lemma}
If $u, v\in\F$ then
\begin{equation}\label{chisum}
\sum_{w\in\F}\chi(w-u)\chi(w-v)=
\begin{cases}-1&{\rm if}\;u\ne v.\\
q-1&{\rm if}\; u=v.
\end{cases}
\end{equation}
\end{lemma}

\noindent{\sl Proof.}
 If $u\ne v$ we have
 \begin{align*}
\sum_{w\in\F}\chi(w-u)\chi(w-v)=&\sum_{w\ne u,v}\chi(w-u)^2\chi(x)
\quad\left(x:=\frac{w-v}{w-u}\right)\\
=&\sum_{w\ne u,v}\chi(x)\\
=&\sum_{x\ne 0,1}\chi(x)\\
=&\sum_{x\ne 0}\chi(x)-\chi(1)\\
=&-1,
\end{align*}
since clearly $\sum_{x\ne 0}\chi(x)=0$. The case $u=v$ is obvious. \hfill$\square$

\medskip

If $q\equiv 1$ mod`$(4)$ then $P(q)$ has adjacency matrix
\[A=\frac{1}{2}(Q-I+J),\]
obtained from $Q$ by replacing every entry $-1$ with $0$. By squaring $A$, we obtain the following lemma:

\begin{lemma}\label{Paleysrg}
 If $q\equiv 1$ mod~$(4)$, and $u$ and $v$ are distinct vertices of $P(q)$, then the number $|(S+u)\cap(S+v)|$ of common neighbours of $u$ and $v$ is
\begin{equation}\label{NoNbrs}
\frac{1}{4}(q-3-2\chi(u-v))=
\begin{cases}
(q-5)/4 & {\rm if }\; u-v\in S,\\
(q-1)/4 & {\rm if }\; u-v\not\in S.
\end{cases}
\end{equation}
\end{lemma}

(When $q\equiv 3$ mod~$(4)$ we find that $|(S+u)\cap(S+v)|=(q-3)/4$ for all pairs $u\ne v$; see also~\cite[Exercise~1.22]{Lem00}.)

This shows that $P(q)$ is a strongly regular graph, with parameters $v=q$ (the number of vertices), $k=(q-1)/2$ (their common valency), $\lambda=(q-5)/4$ and $\mu=(q-1)/4$ (the number of common neighbours of an adjacent or non-adjacent pair of vertices). 

There are several other ways to derive Lemma~\ref{Paleysrg}, for instance from a result of Jacobsthal \cite{Jac06, Jac07} that if $p$ is an odd prime, then the Legendre symbol satisfies
\begin{equation}\label{Jacosum}
\sum_{x=1}^p\left(\frac{x^2+c}{p}\right)=
\begin{cases}
-1 & {\rm if}\; c\not\equiv 0\; {\rm mod}\,(p),\\
p-1 & {\rm if}\; c\equiv 0\; {\rm mod}\,(p);
\end{cases}
\end{equation}
the same proof gives the corresponding result for all finite fields of odd order. One can also deduce Lemma~\ref{Paleysrg} from results of Perron~\cite{Per} and of Kelly~\cite{Kel} on the distribution of quadratic residues (with analogous results for $q\equiv 3$ mod~$(4)$); they both prove their results only in the case where $q$ is prime, though Kelly notes that his arguments are also valid for all odd prime powers, as indeed are those of Perron. Basile and Brutto give a geometric proof in~\cite{BB}.

In fact, by the following lemma one can deduce the strong regularity of $P(q)$, and the values of $k$, $\lambda$ and $\mu$, merely from the facts that $P(q)$ is self-complementary (under the isomorphism $P(q)\to\overline{P(q)},\; v\mapsto av$ for a non-residue $a$), and is arc-transitive (under the automorphisms $v\mapsto av+b$ where $a\in S$ and $b\in\F$).

\begin{lemma}
Any self-complementary arc-transitive graph is strongly regular, with parameters $v=4t+1$, $k=2t$, $\lambda=t-1$ and $\mu=t$ for some integer $t$.
\end{lemma}

\noindent{\sl Proof.} Since the graph is arc-transitive, its automorphism group acts transitively on the vertices, so they all have the same valency $k$. The stabiliser of each vertex is transitive on its neighbours, and hence, since the graph is self-complementary, also on its non-neighbours, so the graph is strongly regular. The complement of a strongly regular graph with parameters $(v, k, \lambda, \mu)$ is also strongly regular, with parameters $(v, v-k-1, v-2-2k+\mu, v-2k+\lambda)$. Since the graph is self-complementary we can equate these parameters, giving
\[v=2k+1\quad{\rm and}\quad \mu=\lambda+1.\]
In any strongly regular graph, double counting of the edges between the neighbours and non-neighbours of a particular vertex gives
\[(v-k-1)\mu=k(k-\lambda-1),\]
so substituting for $v$ and cancelling $k$ gives
\[\mu=k-\lambda-1.\]
Solving the two simultaneous equations for $\lambda$ and $\mu$ give the result, with $t=\mu$. \hfill$\square$

\medskip

In particular, we see that $P(q)$ has parameters $v=q$, $k=(q-1)/2$, $\lambda=(q-5)/4$ and $\mu=(q-1)/4$. 
   

\section{Jacobsthal and sums of squares}

The property of quadratic residues expressed by equation~(\ref{chisum}) can be traced back to the work of Jacobsthal on the representation of primes as sums of squares.

In 1625 Girard stated that each prime $p\equiv 1$ mod~$(4)$ can be written as a sum
\begin{equation}\label{sumsq}
p=a^2+b^2\quad(a, b\in{\mathbb Z})
\end{equation}
of two squares. A few years later Fermat claimed to have a proof, which has never been found, and Euler eventually provided one in 1749. (See~\cite[Ch.~VI]{Dic} for a detailed account of the history of this theorem, and~\cite[\S6.8]{Sti} for a concise summary.) There are several fairly elementary proofs of this result: see~\cite[\S10.1, \S10.6]{JoJo98} or \cite[\S6.7]{Sti}). For instance, one can use simple area calculations to show that if $u^2\equiv -1$ mod~$(p)$ then the lattice $\{(x,y)\in{\mathbb Z}^2\mid y\equiv ux\;{\rm mod}\,(p)\}$ in ${\mathbb R}^2$ has a non-zero element $(a,b)$ within the disc $x^2+y^2<2p$. Zagier has given an elegant one-sentence proof in~\cite{Zag}. However, these proofs are not constructive. There are interesting discussions of constructive proofs in~\cite[\S V.3]{Dav} and~\cite{Wag}; these include Gauss's simple but hardly practical solution
\[a=\bigg\langle\frac{1}{2}{2k\choose k}\bigg\rangle,\quad b=\langle (2k)!\,a\rangle,\]
where $p=4k+1$ and $\langle n\rangle$ is the residue of $n$ mod~$(p)$ closest to $0$, and also in~\cite{Wag} some efficient modern algorithms for solving (\ref{sumsq}).

In 1907 Jacobsthal~\cite{Jac07} published explicit formulae for integers $a$ and $b$ satisfying (\ref{sumsq}), based on work in his thesis~\cite{Jac06}. Specifically, he took 
\[a=\frac{\varphi(r)}{2},\quad b=\frac{\varphi(n)}{2}\]
for any residue $r$ and non-residue $n$ mod~$(p)$, where
\[\varphi(e):=\sum_{m=1}^p\chi(m)\chi(m^2+e).\]
(For typographic convenience we write $\chi$ here, rather than the Legendre symbol used by Jacobsthal.) It is easy to see that $\varphi(e)$ is even. Now $\varphi(e)=0$ if $e\equiv 0$ mod~$(p)$, and otherwise, since $\varphi(e)=\chi(x)\varphi(ex^2)$ for all $x\not\equiv 0$ mod~$(p)$, $\varphi^2(e)$ depends only on whether $e$ is congruent to a residue or a non-residue mod~$(p)$; thus
\[\sum_{e=1}^p\varphi^2(e)=\frac{p-1}{2}\left(\varphi^2(r)+\varphi^2(n)\right).\]
where $r$ and $n$ are any residue and non-residue. On the other hand, calculating the left-hand side directly, using the definition of $\varphi$ and summing first over $e$, leads via equation~(\ref{chisum}) to
\[\sum_{e=1}^p\varphi^2(e)=p(p-1)\left(1+\chi(-1)\right).\]
It follows that if $p\equiv 1$ mod~$(4)$, so that $\chi(-1)=1$, then
\[\left(\frac{\varphi(r)}{2}\right)^2+\left(\frac{\varphi(n)}{2}\right)^2=p,\]
as required.

For example, if $p=13$ we can take $r=1$ and $n=2$, with $\varphi(1)=6$ and $\varphi(2)=-4$, so that $13=3^2+(-2)^2$.

It is tempting to speculate whether Gauss was aware of (\ref{chisum}) in some form, when $q=p$. The author has no direct evidence for this (it does not appear in the {\sl Disquisitiones}), but the Gauss expert Franz Lemmermeyer has commented~\cite{Lem16} that he ``would have seen the proof in a second",  either by counting $\F_p$-rational points on the curve $x^2-y^2=1$, or along the following lines.

Use Euler's criterion $\chi(x)\equiv x^{(p-1)/2}$ mod~$(p)$ and the Binomial Theorem to express each summand in the left-hand side $L$ of~(\ref{chisum}) as a polynomial in ${\mathbb Z}_p[w]$, and then use the fact~\cite[\S 19]{Gau} that if $k\in\mathbb N$ then
\[\sum_{w\in{\mathbb Z}_p}w^k=
\begin{cases}
p-1 & {\rm if} \; k\equiv 0 \; {\rm mod} \, (p-1),\\
0 & {\rm otherwise},
\end{cases}
\]
to show that $L$, regarded as an element of $\mathbb Z$, is congruent to $-1$ mod~$(p)$. If $u\ne v$ then two of the summands in $L$ are equal to $0$, and the remaining $p-2$ are each $\pm 1$, so $|L|\le p-2$ and hence $L=-1$. The case $u=v$ is trivial.


\section{Perron, Brauer, Hopf and Schur}

In 1952 Perron~\cite{Per} studied the distribution of quadratic residues modulo a prime $p$. His theorems that are relevant here can be stated concisely as follows:

\begin{thm}\label{Perronthm}
 Let $\F=\F_p$ for a prime $p=4n\pm 1$, define $S_0=S\cup\{0\}$, and let $a\in\F^*$. Then
\[|(S_0+a)\cap S_0|=
\begin{cases}
n & {\rm if}\; p=4n-1\;{\rm or}\; a\not\in S,\\
n+1 & {\rm if}\; p=4n+1\;{\rm and}\; a\in S.
\end{cases}
\]
\end{thm}
His similar results for the set of non-residues follow immediately on taking complements,  and the corresponding result for $S$ rather than $S_0$ can be deduced from the equation
\[|(S+a)\cap S|=|S_0+a)\cap S_0|-|\{\pm a\}\cap S|.\]

In the same year A.~Brauer~\cite{Bra}, writing on the distribution of quadratic residues with applications to Hadamard matrices (or Hadamard determinants as he called them --- see Section~\ref{Hadamard}), noted that some of Perron's results were corollaries to a theorem in Jacobsthal's thesis~\cite{Jac06}, that if $p$ is an odd prime, and $c\not\equiv 0$ mod~$(p)$, then
\begin{equation}\label{Jacobeqn}
\sum_{x=1}^p\chi(x^2+c))=-1.
\end{equation}
Referring to the case $p=4n-1$ of Theorem~\ref{Perronthm}, Brauer wrote: ``As early as 1920, H.~Hopf showed me his proof of this theorem using (\ref{Jacobeqn}), and its application to the construction of Hadamard determinants of order $p+1=4n$. However he never published it since I.~Schur already knew this result at that time. Independently this theorem and its application to Hadamard determinants were published by R.~E.~A.~Paley~\cite{Pal} in 1933.''

This story was repeated by Dembowski in~\cite[p.~97]{Dem}, where he defined Paley designs as examples of Hadamard designs (see his footnote~(4)). Dembowski asserted that, according to Brauer, Schur already knew of these designs; however, this is not clear, since Brauer did not mention designs in~\cite{Bra}. Certainly the designs are implicit in the matrices, but the first to make an explicit connection seems to have been Todd~\cite{Tod}, in 1933 (see Section~\ref{Todd}).

Since Schur and Frobenius were Jacobsthal's advisors for his 1906 doctoral thesis, one can assume that Schur actually knew the result Brauer refers to much earlier than 1920.
Theorem~\ref{Perronthm} extends in the obvious way to finite fields $\F_q$ of any odd order $q$. In particular, if $q\equiv 1$ mod~$(4)$ we see that
\[|(S+a)\cap S|=
\begin{cases}
(q-5)/4 & {\rm if}\; a\in S,\\
(q-1)/4 & {\rm if}\; a\not\in S,
\end{cases},\]
giving the parameters $\lambda$ and $\mu$ for the strongly regular graph $P(q)$.


\section{Hadamard matrices and designs}\label{Hadamard}

Here we briefly discuss Hadamard matrices and designs, mentioned in the preceding section.  A {\sl Hadamard matrix\/} of order $m$ is an $m\times m$ matrix $H$, with all entries equal to $\pm 1$, and with mutually orthogonal rows, so that $HH^T=mI$. These matrices are named after Jacques Hadamard (1865--1963), who proved in 1893~\cite{Had}  that an $m\times m$ complex matrix $H=(h_{ij})$, with $|h_{ij}|\le 1$ for all $i$ and $j$, satisfies $|\det H|\le  m^{m/2}$; in the case where each $h_{ij}$ is real, $H$ attains this bound if and only if it is  Hadamard matrix. Hadamard matrices have many modern applications, in areas such as engineering, coding theory, cryptography, physics and statistics (see~\cite{HW, Hor, SWW}, for example).

As early as 1867 Sylvester~\cite{Syl}, as part of a study of orthogonal matrices, gave a recursive construction of what later became known as Hadamard matrices, of each order $m=2^e$: he started with $H=I_1=(1)$ and used a Hadamard matrix $H$ of order $m$ to construct the Kronecker product
\[\left(\,\begin{matrix}+&+\cr +&-\cr \end{matrix}\,\right)\otimes H
=\left(\,\begin{matrix}H&H\cr H&-H\cr \end{matrix}\,\right),\]
a Hadamard matrix of order $2m$. (It is typographically convenient to write entries $1$ and $-1$ as $+$ and $-$.) Connoisseurs of Victorian English literary style and social attitudes will, no doubt, appreciate the way in which Sylvester commended his ideas to his readers: he described the many possible applications of his theory, listed at some length in the title of his paper, as ``... furnishing interesting food for thought, or a substitute for the want of it, alike to the analyst at his desk and the fine lady in her boudoir."

In 1898 Scarpis~\cite{Sca} gave a construction of Hadamard matrices of certain orders, but there seems to have been little further progress in their construction until Paley's paper~\cite{Pal} in 1933. This was motivated by problems in combinatorics and geometry, raised by his Cambridge contemporaries Todd~\cite{Tod} and Coxeter~\cite{Cox33} in papers which appeared in the same volume.

It is easy to show~\cite[Lemma~6.28]{JoJo00} that, apart from trivial examples of order $1$ or $2$, if a Hadamard matrix of order $m$ exists then $m\equiv 0$ mod~$(4)$. The Hadamard Conjecture is that the converse is also true. For example, in 1933 Paley wrote~\cite[p.~312]{Pal} ``It seems probable that, whenever $m$ is divisible by $4$,  it is possible to construct an orthogonal matrix of order $m$ composed of $\pm 1$, but the general theorem has every appearance of difficulty.'' (At that time, a square matrix was called orthogonal if it had mutually orthogonal rows; nowadays we impose the extra requirement that all rows have unit length, though it might be more consistent to call such a matrix orthonormal.)

The following process shows that Hadamard matrices lead naturally to certain block designs. Multiplying various rows and columns of a Hadamard matrix by $-1$ yields a {\sl normalised\/} Hadamard matrix, one in which all entries in the first row and column are equal to $1$. Deleting this row and column leaves a square matrix of order $m-1$; its columns and rows can be identified with the points and blocks of a block design, with entries $\pm 1$ indicating incidence or non-incidence of points and blocks. If $m\ge 4$ there are $m-1$ points and blocks, each block has size $(m-2)/2$, and any two blocks have $(m-4)/4$ points in common. A block design with these properties is called a {\sl Hadamard design}. As noted by Todd~\cite{Tod} (see Section~\ref{Todd}), this process is reversible, so that every Hadamard design corresponds to a normalised Hadamard matrix.


\section{Paley}

 

Raymond Edward Alan Christopher Paley was born in Bournemouth, UK, on 7 January 1907, the son of an army officer who died before Paley was born. He attended Eton College, where he was a King's Scholar, entitling him to reduced fees. This school, founded in 1440 by King Henry VI, is noted for having educated 19 British prime ministers, together with one Fields medallist and a number of fictional characters ranging from Captain Hook, via Bertie Wooster, to James Bond. He then studied mathematics at Trinity College, Cambridge, taking his PhD under the supervision of J.~E.~Littlewood. 

In his short life, Paley's main mathematical contributions were in analysis, and many were of considerable significance: Littlewood-Paley theory, the Paley-Wiener Theorem, and the Paley-Zygmund inequality, for example.  A footnote in~\cite[p.~318]{Pal}, citing papers by Walsh, Kaczmarz and himself, suggests that Paley's expertise in constructing Hadamard matrices arose partly  from his work on orthogonal functions. In addition to Littlewood, he collaborated with Zygmund, who spent the year 1930--31 in Cambridge; Zygmund's 1935 book {\sl Trigonometric Series} drew heavily on their joint work. In 1932 Paley obtained a research fellowship to allow him to work with Wiener at MIT. While in the USA, he also collaborated with P\'olya, who was visiting Princeton. Some of his collaboration with Coxeter and Todd, which gave rise to his constructions of Hadamard matrices in~\cite{Pal}, may have taken place in the USA, as they visited Princeton in 1932--33 and 1933--34 respectively.

In the foreword to his edition of {\sl Littlewood's Miscellany}, B\'ela Bollob\'as has written that ``Paley ... was one of the greatest stars in pure mathematics in Britain, whose young genius frightened even Hardy." However, after a highly promising start to his career, Paley died on 7 April 1933 at the age of 26, caught in an avalanche while skiing near Banff.  Wiener wrote in~\cite{Wien} that `` ... he was already recognised as the ablest of the group of young English mathematicians who have been inspired by the genius of G.~H.~Hardy and J.~E.~Littlewood. In a group notable for its brilliant technique, no one had developed this technique to a higher degree than Paley. Nevertheless he should not be thought of primarily as a technician, for with his ability he combined creative power of the first order." MathSciNet lists 23 publications by Paley, including a reprint and a Russian translation of his work with Wiener on Fourier transforms.


\section{Paley's Hadamard matrix constructions}\label{Paleymatrices}

In~\cite{Pal} Paley gave several constructions of Hadamard matrices, based on the combinatorial properties of quadratic residues such as equation~(\ref{chisum}). Perhaps surprisingly, the name Hadamard does not appear in his paper, nor in the accompanying papers by Todd~\cite{Tod} and Coxeter~\cite{Cox33}, apart from once in the title of a bibliographic reference (to an abstract~\cite{Gil} of a talk by Gilman) added by Paley after submitting his paper: Paley and Coxeter referred to what we now call Hadamard matrices as `U-matrices', while Todd had no special name for them.

Paley described two constructions, based on finite fields, which give Hadamard matrices of order $m=q+1$ or $2(q+1)$ for each prime power $q\equiv 3$ or $1$ mod~$(4)$ respectively. He gave proofs only in the cases where $q$ is prime, crediting Todd and Coxeter for the proof when $q\equiv 3$  mod~$(4)$, and Davenport (another Cambridge contemporary) for pointing out the crucial property (\ref{chisum}) of the Legendre symbol; however, later in his paper he noted~\cite[p.~316]{Pal} that his proofs generalise easily to odd prime powers. He then showed in~\cite[Table 1]{Pal} that combinations of these constructions and that of Sylvester yield Hadamard matrices of all orders $m\equiv 0$ mod~$(4)$ up to and including $200$, with the exceptions of $92$, $116$, $156$, $184$ and $188$. Subsequently these and many other orders $m$ have been dealt with, but the conjecture is still open. (Around the same time as Paley, Gilman showed how to construct Hadamard matrices of order $2^{\nu}n_1^{\nu_1}\ldots n_k^{\nu_k}$ where each $n_i\equiv 0$ mod~$(4)$, each $n_i-1$ is prime, and each $\nu_i\ge 1$; reference~\cite{Gil}, cited by Paley, is an abstract of a lecture on this subject given by Gilman, but his work does not seem to have been published.)

Paley stated the above results as lemmas. In another construction, stated as the only theorem in his paper, Paley proved that if $m$ is a power of $2$ one can partition the $2^m$ possible rows of $m$ entries $\pm 1$ into $2^m/m$ sets, each set forming the rows  of a Hadamard matrix of order $m$ (see Section~\ref{Paley3} for an outline proof of this result).


\subsection{Paley's first construction}\label{Paley1}

Let $Q$ be the Jacobsthal matrix for the field $\F=\F_q$, where $q$ is an odd prime power, and let $R$ be the row vector $(1,1,\ldots,1)$ of length $q$. In~\cite[Lemma~2 and p.~316]{Pal}, Paley showed that if $q\equiv 3$ mod~$(4)$ then
\[H=\left(\,\begin{matrix}1&R\cr R^T&Q-I\cr \end{matrix}\,\right)\]
is a Hadamard matrix of order $m=q+1$, where $I$ denotes the identity matrix of order $q$. The fact that distinct rows are orthogonal follows immediately from equation~(\ref{chisum}). These matrices are now known as {\sl Paley-Hadamard matrices of type I}.


\subsection{Paley's second construction}\label{Paley2}

In~\cite[Lemma~3]{Pal} Paley started with the matrix
\[(B_{ij})=\left(\,\begin{matrix}0&R\cr R^T&Q\cr \end{matrix}\,\right),\]
where $q\equiv 1$ mod~$(4)$, and then replaced each entry $B_{ij}=\pm 1$ or $0$ with the $2\times 2$ matrix
\[\pm\left(\,\begin{matrix}+&+\cr +&-\cr \end{matrix}\,\right)\quad{\rm or}
\quad \left(\,\begin{matrix}+&-\cr -&-\cr \end{matrix}\,\right)\]
respectively.  This gives a symmetric Hadamard matrix of order $m=2(q+1)$. The proof is similar to that for his first construction. These matrices are now known as {\sl Paley-Hadamard matrices of type II}.


\subsection{Paley's third construction}\label{Paley3}

In another construction, stated as the only theorem in his paper, Paley proved that if $m$ is a power of $2$ one can partition the $2^m$ possible rows of $m$ entries $\pm 1$ into $2^m/m$ sets, each set forming the rows  of a Hadamard matrix of order $m$. Since this result seems to be rather less well-known, and since it was subsequently used by Todd~\cite{Tod} and Coxeter~\cite{Cox33}, we give an outline proof, following Paley's notation (though with more modern terminology).
 
Let $m=2^k$, let $i, j\in\{0, 1, \ldots, 2^k-1\}$ have binary representations
\[i=\sum_{\lambda=0}^{k-1}\eta_{\lambda}2^{\lambda},\quad j=\sum_{\lambda=0}^{k-1}\zeta_{\lambda}2^{\lambda}\]
where $\eta_{\lambda}, \zeta_{\lambda}\in\{0, 1\}$, and write
\[A_{ij}=\prod_{\lambda=0}^{k-1}(-1)^{\eta_{\lambda}\zeta_{\lambda}-1-\lambda}.\]
The matrix $M=(A_{ij})$ is the $k$th Kronecker power of the matrix ${+\; +\choose +\; -}$.

Then $A_{i_1j}A_{i_2j}=A_{ij}$, where the binary representation of $i$ is the term by term mod~$(2)$ sum of those for $i_1$ and $i_2$. Thus the $m$ rows $(A_{ij}),\; i$ fixed, $j=0,\ldots, m-1$, form an elementary abelian group under term by term multiplication. This multiplication rule also shows that distinct rows of $M$ are orthogonal, so $M$ is a Hadamard matrix.

For each of the $2^m$ sequences $B=(B_j)\in\{\pm 1\}^m$ we have a Hadamard matrix $(A_{ij}B_j)$. As shown by Paley, if two of these matrices have a row in common (possibly in different positions), they have all their rows in common, so these matrices $(A_{ij}B_j)$ are partitioned into $2^m/m$ sets of size $m$, those in the same or different sets having all or none of their rows in common. Choosing one matrix from each set proves the theorem.


\subsection{Todd's paper}\label{Todd}

Todd used Paley's theorem in the related paper~\cite{Tod}. Motivated by work of Coxeter~\cite{Cox33} on polytopes, Todd was interested in the problem of finding $4n-1$ subsets of size $2n-1$ in a $(4n-1)$-element set, each pair having an intersection of size $n-1$. (Such an arrangement is now called a Hadamard design with parameters $(4n-1, 2n-1, n-1)$, and the chosen subsets are called blocks.) As explained in Section~\ref{Hadamard}, this is equivalent to finding a Hadamard matrix $H$ of order $m=4n$. Normalising $H$, then removing the first row and column, and finally replacing all entries $-1$ with $0$ gives the incidence matrix of the design.

One obvious solution to Todd's problem, corresponding to the matrix $M$ of order $m=2^k$ in the proof of Paley's main theorem~\cite{Pal} (see Section~\ref{Paley3}), is to take the blocks to be the hyperplanes in the $(k-1)$-dimensional projective geometry $PG(k-1,2)$ over the field $\F_2$, with $n=2^{k-2}$; the points and lines of the Fano plane correspond to the simplest case, $k=3$. Another solution, corresponding to Paley's Lemmas~2 and 4, is to take the blocks to be the translates of the set $S$ of quadratic residues in the field $\F_q$ of order $q\equiv 3$ mod~$(4)$, with $n=(q+1)/4$. (Dembowski called this a {\sl Paley design}~\cite[p.~97]{Dem}.) Todd also constructed other examples, for instance classifying them all in the case where there are $4n-1=15$ points.

In addition, Todd considered the automorphism group of such a design, that is, the largest subgroup of $S_{4n-1}$, acting on the points, which permutes the blocks. In the finite geometry example, this is the collineation group $PGL_k(2)$ of the geometry (see~\cite[p.~31]{Dem} for notation for groups of projective transformations). In the quadratic residue example, the automorphism group contains the subgroup
\[A\Delta L_1(q):=\{v\mapsto av^{\gamma}+b\mid a\in S,\, b\in\F_q,\, \gamma\in{\rm Gal}\,\F_q\}\]
of index $2$ in $A\Gamma L_1(q)$ (see also Section~\ref{charautos} for this group); in some cases,  such as when $q=19$, $23$ or $27$ (in fact for all $q\ge 19$, by a later result of Kantor~\cite{Kan69}) this is the whole automorphism group, but in other cases the automorphism group is larger, for instance isomorphic to $PSL_2(q)$ (acting with degree $q$) when $q=7$ or $11$. (Although Todd did not mention this, when $4n-1=7$ the isomorphism between the projective geometry and quadratic residue solutions illustrates the isomorphism $PGL_3(2)\cong PSL_2(7)$.) 


\subsection{Coxeter's paper}\label{Coxeter}

In the other related paper~\cite{Cox33}, Coxeter was interested in generalising the well-known partitions of the vertices of the cube or of the dodecahedron into those of two or five tetrahedra. Schoute~\cite{Sch} had already given similar examples of such compound polytopes in dimension $4$, and Coxeter wanted to construct examples in dimensions $m\ge 5$. In this case the possibilities are more restricted, since the only regular polytopes are the simplex $\alpha_m$, the cross-polytope $\beta_m$, and its dual,  the `measure polytope', or $m$-cube $\gamma_m$. (When $m\ge 5$ there are no analogues of the dodecahedron $\{5, 3\}$ and the icosahedron $\{3,5\}$ for $m=3$, or of the 24-cell $\{3,4,3\}$, the 120-cell $\{5,3,3\}$ and the 600-cell $\{3,3,5\}$ for $m=4$; here the brackets $\{\ldots\}$ denote Schl\"afli symbols for polytopes, see~\cite{Cox73}.)

\begin{figure}[h!]
\begin{center}
\begin{tikzpicture}[scale=0.2, inner sep=0.8mm]

\node (A) at (15,9) {};
\node (B) at (6,15) {};
\node (C) at (-15,11) {};
\node (D) at (-6,5) {};
\draw [thick, red] (15,9) to (6,15) to (-15,11) to (-6,5) to (15,9);

\draw [thick, red] (-15,-9) to (-6,-15) to (15,-11);
\draw [thick, red] (15,9) to (15,-11) {};
\draw [thick, red] (-15,11) to (-15,-9) {};
\draw [thick, red] (-6,5) to (-6,-15) {};

\draw [thick] (15,9) to (-15,11) to (-6,-15) to (15,9);
\draw [thick] (-15,-9) to (-6,5) to (15,-11);
\draw [thick] (6,15) to (-6,5);
\draw [thick] (-10.5,-2) to (4.5,-3) to (0,10);
\draw [thick] (6,15) to (7.9,9.5);
\draw [thick] (10.5,2) to (15,-11);

\draw (15,9) to (10.5,2) to (4.5,-3);
\draw (-6,-15) to (0,-10) to (4.5,-3);
\draw (15,-11) to (0,-10);
\draw (-15,-9) to (-8,-9.5);

\path [fill=lightgray, thin] (15,9) to (0,10) to (1.3,6.4) to (15,9); 
\path [fill=lightgray, thin] (15,9) to (1.3,6.3) to (4.5,-3) to (15,9); 

\draw [fill=lightgray, thin] (15,9) to (10.5,2) to (4.5,-3) to (15,9);

\path [fill=lightgray, thin] (4.5,-3.1) to (-5.9,-15) to (0,-9.9) to (4.5,-3.1);
\path [fill=lightgray, thin] (4.4,-3) to (-5.9,-14.9) to (-5.9,-2.25) to (4.4,-3);

\path [fill=lightgray, thin] (-6.1,-15) to (-6.1,-2.25) to (-10.5,-2) to (-6,-15);
\path [fill=lightgray, thin] (-14.9,11) to (-5.9,5.1) to (-0.1,10) to (-14.9,11);
\path [fill=lightgray, thin] (-15,10.9) to (-6.1,4.9) to (-10.5,-1.9) to (-15,10.9);

\end{tikzpicture}
\end{center}
\caption{A stella octangula, inscribed in a cube} 
\label{StellaOct}
\end{figure}
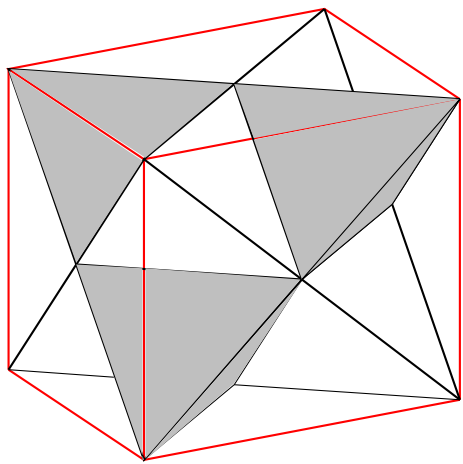

Coxeter defined a {\sl compound polytope\/} in dimension $m$ to be a set of $D$ concentric, finite, convex, $m$-dimensional polytopes $\Pi$, which are transitively permuted by the symmetry group of the set. Let $\Pi$ have $V$ vertices and $F$ faces. The compound polytope is {\sl vertex-regular\/} if the $DV$ vertices of its components are the $v$ vertices of a regular polytope $\pi_1$, each taken $d_1$ times for some $d_1\ge 1$, so that $DV=d_1v$. For example, when $m=3$ one could take $\pi_1$ to be a cube $\gamma_3=\{4,3\}$ (see Figure~\ref{StellaOct}) or dodecahedron $\{5,3\}$, with its vertices partitioned into those of $D=2$ or $5$ tetrahedra $\Pi=\alpha_3=\{3,3\}$; here $V=F=4$, $v=8$ or $20$, and $d_1=1$. For an example where $d_1>1$ take the two possible partitions for the dodecahedron, mirror-images of each other, so that $D=10$ and $d_1=2$.

Dually, Coxeter defined a compound to be {\sl face-regular\/} if the $DF$ bounding spaces (hyperplanes spanned by faces) of its components are the $f$ bounding spaces of a regular polytope $\pi$, each taken $d$ times, so that $DF=df$. For example, $\pi$ could be the octahedron $\pi=\beta_3=\{3,4\}$ or icosahedron $\{3,5\}$, regarded as the intersection of $D=2$ or $5$ tetrahedra, so that $d=1$. In the latter case one could also take the two possible sets of five tetrahedra, giving $D=10$ and $d=2$.

A simple argument, counting vertices, shows that if $m\ge 5$ then in any vertex-transitive compound, $\pi_1$ must be $\gamma_m$, with a dual result $\pi=\beta_m$ for face-transitive compounds. It follows that the only possibilities are (using a slightly more consistent version of Coxeter's notation):
\begin{enumerate}
\item $d_1\gamma_m[D\beta_m]$, a vertex-transitive compound of $D$ cross-polytopes $\Pi=\beta_m$, their vertices $d_1$ at a time forming an $m$-cube $\pi_1=\gamma_m$;
\item its dual $[D\gamma_m]d\beta_m$, a face-transitive compound of $D$ hypercubes $\Pi=\gamma_m$, their bounding spaces $d$ at a time forming a cross-polytope $\pi=\gamma_m$;
\item $d_1\gamma_m[D\alpha_m]d\beta_m$, a self-dual vertex- and face-transitive compound of $D$ simplices $\Pi=\alpha_m$, their vertices and faces $d_1$ and $d$ at a time forming $\pi_1=\gamma_m$ and $\pi=\beta_m$.
\end{enumerate}
For example, Kepler's stella octangula, shown in Figure~\ref{StellaOct}, is the compound
\[\gamma_3[2\alpha_3]\beta_3=\{4,3\}[2\{3,3\}]\{3,4\}.\]
This is formed as above from two tetrahedra $\alpha_3=\{3,3\}$, inscribed in a cube $\gamma_3=\{4,3\}$ and intersecting in an octahedron $\beta_3=\{3,4\}$, with $d_1=d=1$.

Coxeter showed that compounds of type (1) and (2) exist, in dimension $m$, if and only if a compound of type (3) exists, in dimension $m-1$, with the same values of $d_1$, $d$ and $D$. Moreover, the existence of such compound polytopes is equivalent to that of certain Hadamard matrices, as follows.

The points in ${\mathbb R}^m$ with all coordinates $\pm 1$ form the vertices of an $m$-cube $\gamma_m$. If $H$ is a Hadamard matrix of order $m$, then the rows of $H$ and $-H$ are the vertices of a cross-polytope $\beta_m$ inscribed in $\gamma_m$, and this is just one component of a compound polytope of type (1); the dual compound has type (2). Similarly, given a normalised Hadamard matrix of order $m$, deleting the first column gives $m$ rows, the vertices of a simplex $\alpha_{m-1}$ inscribed in $\gamma_{m-1}$, and this leads to a compound of type (3), but in dimension $m-1$. The converse is also true, that for each of these three types, any such compound arises in this way from a Hadamard matrix.

Thus, apart from trivial cases where $m\le 2$, $m$ must be divisible by $4$. Conversely, Coxeter deduced from Paley's work~\cite{Pal} on Hadamard matrices that compounds of types~(1) and (2) exist in dimensions $4, 8, 12, \ldots, 88$, and those of type~(3) exist in dimensions $3, 7, 11, \ldots, 87$, although the next cases $92$ and $91$ were at that time still undecided. (A Hadamard matrix of order $92$ was found in 1962 by Baumert, Golomb and Hall~\cite{BGH}, using a construction due to Williamson~\cite{Wil} and a great deal of computing.) Similarly,  when $m=2^k$ the $2^m/m$ Hadamard matrices given by Paley's main theorem yield a compound
\[\gamma_{2^k-1}[2^{2^k-k-1}\alpha_{2^k-1}]\beta_{2^k-1}\]
for each $k\ge 2$, generalising the stella octangula for $k=2$.

As shown by Todd~\cite{Tod}, Hadamard matrices of order $m=4n$ are equivalent to certain block designs, or arrangements of subsets, on $4n-1$ points, so Todd's results on the orders of their automorphism groups can be used to consider the possibilities for the parameters $d_1$ and $D$. If the automorphism group has order $N$ there are $(4n-1)!/N$ different designs. If we take the $4n-1$ elements as coordinate places, each block corresponds to a point in ${\mathbb R}^{4n-1}$ with coordinates respectively $\pm 1$ at its elements and non-elements. Each design therefore gives $4n-1$ points which, together with $(1,1,\ldots,1)$, are the vertices of a simplex $\alpha_{4n-1}$ inscribed in $\gamma_{4n-1}$. The $(4n-1)!/N$ different designs thus yield  a compound of each type (1), (2) and (3), with $d_1=(4n-1)!/N$ and hence
\[D=\frac{d_1v}{V}=\frac{2^m\,d_1}{2m}=\frac{2^{4n-3}\,d_1}{n}=\frac{2^{4n-3}\,(4n-1)!}{nN}.\]
For example, for $n=2$ the design based on quadratic residues mod~$(7)$ has automorphism group $PSL_2(7)$ of order $N=168$, giving such compounds with $d_1=30$ and $D=480$. Similarly, for $n=3$, using quadratic residues mod~$(11)$ gives $N=|PSL_2(11)|=660$, $d_1=60480$ and $D=10321920$.

In his text~\cite{Cox73} on regular polytopes, Coxeter considered compound polytopes of dimensions $3$ and $4$ in considerable detail, but unfortunately when it came to higher dimensions~\cite[pp.~287--288]{Cox73} he wrote simply ``To save space, we have disregarded the possibility of compounds in more than four dimensions", and after giving a few examples, ``The theory of these compounds is connected with orthogonal matrices of $\pm 1$'s", with no references to Hadamard, Sylvester, Paley or Todd.


\section{The origin of the Paley graphs}

The concept of a graph does not appear in Paley's paper~\cite{Pal}, nor in the accompanying papers by Todd~\cite{Tod} and Coxeter~\cite{Cox33}. This is not surprising, since they had no need of graphs for their work, and in any case graphs were little known and rarely studied in the 1930s; indeed, the first textbook on graph theory, by K\"onig~\cite{Kon}, was not published until 1936, three years after Paley's death. The graphs which eventually carried his name first appeared in the literature nearly 30 years later, in two highly influential and almost simultaneous papers, one by Sachs~\cite{Sac}, and the other by Erd\H os and R\'enyi~\cite{ER}.

\subsection{Sachs}\label{Sachs}

In 1962 Sachs~\cite{Sac} introduced the concept of a self-complementary graph, one which is isomorphic to its complement. In this paper he was particularly interested in such graphs which are also regular (meaning that all vertices have the same valency), and cyclic (invariant under a cyclic permutation of the vertices). As examples, he constructed the Paley graphs $P(q)$ for primes $q=p\equiv 1$ mod~$(4)$, using elementary properties of quadratic residues and Legendre symbols to verify that they satisfy these conditions. For instance, multiplying all vertices by a fixed non-residue induces an isomorphism $P(q)\to\overline{P(q)}$. Sachs did not consider the automorphism group, apart from noting that translation by $1$ confirms the cyclic property. Nor did he consider the general case of prime powers $q$ (in this case, $P(q)$ is again self-complementary and regular, but it is not cyclic unless $q$ is prime). He did not name these graphs in his paper, nor did he cite Paley, or any other source, for them. 


\subsection{Erd\H os and R\'enyi }

In 1963 Erd\H os and R\'enyi~\cite{ER} considered the following problem (among many others): given an asymmetric graph $G$ (one with no non-identity automorphisms), what is the minimum number $A(G)$ of edge-changes (deletions or insertions of edges) required to allow a non-identity automorphism? If $\Delta_{uv}$ denotes the number of vertices $w\ne u, v$ adjacent to just one of the distinct vertices $u$ and $v$ in $G$, then deleting those vertices $w$ allows an automorphism transposing $u$ and $v$, and fixing all other vertices, so
\begin{equation}\label{AleDelta}
A(G)\le \min_{u\ne v}\Delta_{uv}.
\end{equation}
Simple counting arguments show that if $G$ has order $n$ then
\begin{equation}\label{Deltabd}
\min_{u\ne v}\Delta_{uv}\le\lfloor\frac{n-1}{2}\rfloor,
\end{equation}
so that
\begin{equation}\label{Abd}
A(G)\le \lfloor\frac{n-1}{2}\rfloor.
\end{equation}
(See~\cite{ABMSST} for a much stronger bound $A(G)\le 5$ for planar graphs.) In their paper, Erd\H os and R\'enyi conjectured 
that no asymmetric graph attains the upper bound in~(\ref{Abd}). However, they constructed the Paley graphs $P(q)$ (without referring to Paley) as examples (far from asymmetric) of what they called $\Delta$-graphs, those which attain equality in~(\ref{Deltabd}): indeed, it follows immediately from~(\ref{NoNbrs}) that $\Delta_{uv}=(q-1)/2$ for all pairs $u\ne v$ in $P(q)$. 

They constructed $P(q)$ first~\cite[p.~301]{ER} for primes $q\equiv 1$ mod~$(4)$, referring to Lagrange, Perron~\cite{Per} and Kelly~\cite{Kel} for the fact, equivalent to~(\ref{NoNbrs}), that if $a\ne 0$ in $\F$ then $S+a$ contains $(q-1)/4$ quadratic non-residues. (It is frustrating that they gave no citation for Lagrange.) 
Later~\cite[p.~302]{ER} they extended the construction to all prime powers $q\equiv 1$ mod~$(4)$, again quoting Kelly~\cite{Kel} for the required form of~(\ref{NoNbrs}) in this more general context. They did not consider the full automorphism group of $P(q)$, merely noting that translation by $1$ is always a non-identity automorphism.

In remarks added after submission, Erd\H os and R\'enyi referred to a forthcoming paper by Bose (presumably~\cite{Bos}) on strongly regular graphs, and pointed out that any $\Delta$-graph is strongly regular, of order $n\equiv 1$ mod~$(4)$ and valency $k=(n-1)/2$, though their assertion about the numbers $\lambda$ and $\mu$ (in modern notation) of common neighbours of two adjacent and non-adjacent vertices is clearly incorrect (see~(\ref{NoNbrs}) for the correct values). They also remarked that in the case where $q$ is prime, their $\Delta$-graphs $P(q)$ coincide with the self-complementary graphs constructed by Sachs in~\cite{Sac} (see Section~\ref{Sachs}).


\subsection{Naming the Paley graphs}

Although the graphs $P(q)$ were known, and in the literature, by the early 1960s, Paley's name does not seem to have been associated with them until the early 1970s. The first appearance in the literature the author has found for the term `Paley graph' is in the book by Cameron and van Lint~\cite[p.~14]{CvL75}, published in 1975, where it is introduced as if it was already accepted terminology. It appears in the 2nd edition (1993) of the book by Biggs~\cite{Big} on Algebraic Graph Theory, but not the first, published in 1974. The earliest title in MathSciNet containing the term is a paper~\cite{BEH} by Blass, Exoo and Harary, published in 1981. On the other hand, the term `Paley design' was used by Dembowski in 1968, in his book on Finite Geometries~\cite{Dem}, and the following year by Kantor in~\cite{Kan69}, so perhaps the term `Paley graph' evolved naturally from this. Several colleagues have suggested Jaap Seidel as the originator of the term, but he is unfortunately no longer with us to confirm or deny this. His highly influential 1965 paper with van Lint~\cite{vLS} cited that of Paley~\cite{Pal}, and used Paley's matrices to construct equilateral point sets in elliptic geometry, while Andries Brouwer has confirmed that the term was standard and understood by all in Eindhoven in the 1970s. Whoever originated the term, whether or not it is justified is discussed in Section~\ref{Attrib}.


\subsection{Pseudo-Paley graphs}

The Paley graphs are strongly regular graphs with parameters $v=q$, $k=(q-1)/2$, $\lambda=(q-5)/4$, $\mu=(q-1)/4$. However, these properties do not characterise them. A {\sl pseudo-Paley graph\/} is a strongly regular graph with the same parameters $v, k, \lambda, \mu$ as a Paley graph. In 2001 Peisert~\cite{Pei} constructed a new infinite class of such graphs, now called {\sl Peisert graphs}, as follows.

Let $\F=\F_q$ where $q$ is an even power of a prime $p\equiv 3$ mod~$(4)$, and let $P^*(q)$ be the Cayley graph for the additive group of $\F$ with respect to the generating set
$\{\omega^j\mid j\equiv 0\;\hbox{or}\;1\; {\rm mod}\, (4)\}$, where $\omega$ is a primitive root in $\F$ (that is, a generator of the group $\F^*$). This is an undirected graph, which is (up to isomorphism) independent of the choice of $\omega$. It is a pseudo-Paley graph, and like $P(q)$ it is self-complementary and arc-transitive, but it is not isomorphic to a Paley graph. Indeed, Peisert showed that, apart from one other graph of order $23^2$, the graphs $P(q)$ and $P^*(q)$ are the only pseudo-Paley graphs which are self-complementary and arc-transitive.


More recently Klin, Kriger and Woldar~\cite{KKW} have used association schemes based on affine planes to construct pseudo-Paley graphs of order $q=p^2$ for odd primes $p$. Most of these are neither self-complementary (for $p\ge 17$) nor arc-transitive (for $p\ge 11$). The numbers of both self-complementary and non-self-complementary examples grow rapidly as $p\to\infty$.


\section{The automorphism group of a Paley graph}

\subsection{Characterising the automorphisms}\label{charautos}

When they introduced the graph $P(q)$, both Sachs~\cite{Sac} and Erd\H os and R\'enyi \cite{ER} noted that $v\mapsto v+1$ is an automorphism of the graph. In Sachs's case $q$ is prime, so this implies that the graph is vertex-transitive, but neither of these papers contains any further discussion of the automorphism group.

In fact, it is clear from its construction that $P(q)$ is invariant under translation by any element of $\F$, under multiplication by any element of $S$, and under any field automorphism of $\F$. For any odd prime power $q$, these transformations generate the subgroup
\[A\Delta L_1(q):=\{v\mapsto av^{\gamma}+b\mid a\in S,\, b\in\F,\, \gamma\in{\rm Gal}\,\F\}\]
of order $q(q-1)e/2$ and of index $2$ in $A\Gamma L_1(q)$, so for each $q\equiv 1$ mod~$(4)$ we have
\[A\Delta L_1(q)\le{\rm Aut}\,P(q).\]

For example, the automorphisms of $P(9)$ induced by the additive, multiplicative and Galois groups of the field $\F=\F_9$ can be seen in Figure~\ref{P(9)} as translations, rotations about $0$ and reflection in the horizontal axis. Similarly, automorphisms $v\mapsto 4v$ and $v\mapsto v+1$ of order $6$ and $13$ of $P(13)$, generating $A\Delta L_1(13)$, can be seen in Figures~\ref{P(13)torus} and \ref{P(13)}.

In fact, the elements of $A\Delta L_1(q)$ are the only automorphisms:

\begin{thm}\label{Paleyautos}
If $q\equiv 1$ {\rm mod}~$(4)$ then
\[{\rm Aut}\,P(q)=A\Delta L_1(q).\]
\end{thm}

\noindent{\sl Proof.} We have already established one inclusion. To prove the reverse inclusion, let $\alpha$ be any automorphism of $P(q)$. Since the subgroup $A\Delta L_1(q)$ of ${\rm Aut}\,P(q)$ acts transitively on the arcs of $P(q)$, by composing $\alpha$ with a suitable element of this subgroup we may assume that $\alpha$ fixes $0$ and $1$. As an automorphism of $P(q)$, $\alpha$ satisfies
\begin{equation}\label{chiinv}
\chi(\alpha(u)-\alpha(v))=\chi(u-v)
\end{equation}
for all $u,v\in\F$. In 1960 Carlitz~\cite{Car} proved that if $q$ is a power of an odd prime $p$, then any permutation of $\F_q$ fixing $0$ and $1$ and satisfying (\ref{chiinv}) has the form $v\mapsto v^{p^i}$ for some $i$. This implies that $\alpha$ is a field automorphism, so $\alpha\in A\Delta L_1(q)$. \hfill$\square$

\medskip

In the terminology introduced by Wielandt in~\cite{Wiel69}, Theorem~\ref{Paleyautos} asserts that the permutation group $A\Delta L_1(q)$ is $2$-closed, that is, it is the full automorphism group of the set of binary relations on $\F_q$ which it preserves.

To be more precise about~\cite{Car}, the theorem Carlitz proved was as follows:

\begin{thm}\label{Carlitzthm}
If $q$ is a power of an odd prime $p$, then any permutation polynomial $\F_q\to\F_q$, which fixes $0$ and $1$ and satisfies~(\ref{chiinv}), has the form $v\mapsto v^{p^i}$ for some $i$.
\end{thm}

However, simple counting shows that any function $\F_q\to\F_q$ can be represented by a polynomial (of degree less than $q$), and as Carlitz later wrote (see his errata~\cite[p.~999]{Car} and Hall's review~\cite{HalMR}), any function satisfying~(\ref{chiinv}) must be a permutation, so his theorem actually applies to any function fixing $0$ and $1$ and satisfying~(\ref{chiinv}). 

Carlitz's theorem was motivated by a problem in finite geometry. In~\cite{Car}, he simply wrote that it ``answers a question raised by W.~A.~Pierce in a letter to the writer," without giving any further details. However Hall, in his review~\cite{HalMR} of the paper, stated that ``The result is somewhat negative in its applications to the theory of projective planes since it shows that a construction by Pierce, analogous to the Moulton construction of non-Desarguesian planes, can yield only Desarguesian planes in the prime case." Clearly this was a reference to Pierce's paper~\cite{Pie}, which was published a year later and which used Carlitz's result to extend Moulton's construction~\cite{Mou}. (For further background, see comments in the introduction of the paper~\cite{McC} by McConnel, who was a student of Carlitz.) 


\subsection{Carlitz}

Leonard Carlitz (1907--1999)  completed his doctorate at the University of Pennsylvania in 1930. After a year working with E.~T.~Bell at Caltech, he spent the academic year 1931--32 as an International Research Fellow in Cambridge, where Hardy had just returned after eleven years in Oxford. According to Hayes's obituary of Carlitz~\cite{Hay},  ``This was the era when Hardy and Littlewood led one of the great centres of research in number theory, and Carlitz found the mathematical atmosphere there exhilarating. His work in additive number theory derives from that period." In~\cite{Lit} Bollob\'as has written ``In December 1931, Hardy and Littlewood announced weekly meetings of a conversation class to start in January 1932 in Littlewood's rooms. According to E.~C.~Titchmarsh, `this was a model of what such a thing should be. Mathematicians of all nationalities and ages were encouraged to hold forth on their own work, and the whole exercise was conducted with a delightful informality that gave ample scope for free discussion after each paper.' Nevertheless, as Dame Mary Cartwright wrote, a little later there was a metamorphosis of Littlewood's conversation class into a larger gathering run by Hardy."

Paley was still in Cambridge that year, so it seems inevitable that, as common members of the group around Hardy and Littlewood, he and Carlitz would have met and come to know each other. To what extent, if any, they influenced each other, is unknown. However, it seems likely that when, nearly 30 years later, Carlitz proved Theorem~\ref{Carlitzthm}, he was completely unaware of any possible connection with Paley and his work.

After his year in Cambridge, Carlitz took up a position at the recently-founded Duke University, in North Carolina. Indeed, it seems likely that it was a strong reference from Hardy which got him this position: who else could have been the ``Oxford don" who, according to Durden~\cite{Dur} (see also~\cite{MacT}), wrote that Carlitz was ``fully master of the technique of his trade" and ``better equipped in the analytic theory of numbers than anyone else in America"? Carlitz spent the rest of his career at Duke University, editing the Duke Mathematical Journal and becoming one of the most prolific mathematicians of the 20th century.

\subsection{Subsequent proofs}

Carlitz's proof of Theorem~\ref{Carlitzthm} is technically quite difficult, and it involves no graph theory or group theory, just calculations with polynomials over finite fields.  However, in the case $q=p$, Theorem~\ref{Paleyautos} is a straightforward consequence of Burnside's theorem~\cite{Bur} that a simply transitive group of prime degree $p$ must be solvable, and hence (by a result of Galois, see~\cite[Satz II.3.6]{Hup}) a subgroup of $AGL_1(p)$, as pointed out by Bruen~\cite{Bru} in 1972. More generally, since ${\rm Aut}\,P(q)$ contains $A\Delta L_1(q)$ and cannot be doubly transitive, it is a rank $3$ permutation group with suborbit-lengths $1, (q-1)/2, (q-1)/2$. It is therefore primitive, since $1+(q-1)/2$ does not divide $q$. In the case $q=p^2$ Theorem~\ref{Paleyautos} therefore follows easily from Wielandt's classification~\cite{Wiel69} of simply primitive groups of degree $p^2$ (see also~\cite[Theorem~B$'$]{JS}): these are either rank $3$ subgroups of the wreath product $S_p\wr S_2$, with suborbit-lengths $1, 2(p-1), (p-1)^2$, or subgroups of $AGL_2(p)$. Comparing suborbit-lengths rules out the first possibility, and in the second case linear algebra gives the result. 
See also results of Dobson and Witte~\cite{DW} on automorphism groups of graphs with $p^2$ vertices. Later we will show how group theory can also deal with higher powers of $p$.

It took some time before the significance of Theorem~\ref{Carlitzthm} for Paley graphs was realised. In 1969 Kantor used it, and cited~\cite{Car}, in proving~\cite[Corollary 8.2]{Kan69}; like Carlitz's theorem, his result was purely about permutations of finite fields, and even though it was in a paper on automorphisms of designs there was no application to automorphism groups. As late as 1972, Shult, having defined $P(q)$ (but not named it) in~\cite[Example~2]{Shu}, wrote that it was an open question whether ${\rm Aut}\,P(q)$ was equal to the group we have called $A\Delta L_1(q)$ or larger, though he noted that certain cases could be handled, citing Higman's paper~\cite{Hig}. Shult gave credit to Kantor for this example, so clearly the link between Carlitz's theorem and the Paley graphs was not widely understood among group theorists in the early 1970s. Dembowski~\cite{Dem}, writing in 1968, included Carlitz's paper~\cite{Car} in his bibliography, possibly to support a citation of Pierce's paper~\cite{Pie} on p.~233, but after two careful searches of the whole book the present author has not found any citation of~\cite{Car}.


\subsection{Generalisations of the main theorem}

In 1963 McConnel~\cite{McC} generalised Theorem~\ref{Carlitzthm}, for any prime power $q$, as follows:

\begin{thm}\label{McCthm}
Let $d$ be a proper divisor of $q-1$, and for $x\in\F$ let $\phi(x):=x^m$ where $m=(q-1)/d$. Then a function $f:\F\to\F$ satisfies
\begin{enumerate}
\item $f(0)=0$ and $f(1)=1$,
\item $\phi(f(u)-f(v))=\phi(u-v)$ for all $u, v\in \F$,
\end{enumerate}
if and only if $f(x)=x^{p^j}$ for some $j$ where $d$ divides $p^j-1$.
\end{thm}

In particular, if $q$ is odd and $d=2$, we have $\phi=\chi$, giving Theorem~\ref{Carlitzthm}. As in~\cite{Car}, neither graph theory nor group theory were used in~\cite{McC}.

In 1972, using the fact that the functions $f$ satisfying condition~(2) of Theorem~\ref{McCthm}, or equivalently satisfying
\[\frac{f(u)-f(v)}{u-v}\in D:=\{x\mid x^m=1\} \; \hbox{for all} \; u\ne v\in \F, \]
form a group under composition, Bruen gave a simple algebraic proof of McConnel's theorem in the case where $q=p$. This was based (as in the case $d=2$) on Burnside's theorem on permutation groups of prime degree. Bruen did not explicitly name or describe this group (let us call it $G(d)$), but it is clear that its elements are the $mqh$ transformations
\[x\mapsto ax^{p^j}+b\quad (a\in D, \; b\in\F ,\;  d\mid p^j-1),\]
where $q=p^e$ and $h=\gcd(m,e)$; those elements also satisfying condition~(1), or equivalently $a=1$ and $b=0$, form the subgroup fixing $0$ and $1$.

In 1973 Bruen and Levinger~\cite{BL} extended Bruen's algebraic proof of McConnel's theorem to the case of all prime powers $q$, using ideas taken from Wielandt's proof of Burnside's theorem given in Passman's book~\cite[Theorem 7.3]{Pas}. (Dress, Klin and Muzychuk have given an elementary and largely geometric proof of Burnside's theorem in~\cite{DKM}, together with a detailed survey of alternative proofs by Burnside, Schur, Wielandt and others.) An essential ingredient in Bruen and Levinger's proof is the vector space of all functions $\F\to\F$ (represented as polynomials of degree less than $q$), and the actions on it of various groups of permutations of $\F$. Some of these ideas overlap with those involving invariant relations and functions, developed by Wielandt in~\cite{Wiel69}. Taking $d=2$, the paper~\cite{BL} seems to be the first to give an explicit description of the elements of ${\rm Aut}\,P(q)$, and thus to give an implicit statement of Theorem~\ref{Paleyautos}.

In 1990 Lenstra~\cite{Len} gave a rather shorter proof of McConnel's theorem, based on that of Bruen and Levinger. He stated the theorem in the slightly more elegant form that
\[G(d)=\{x\mapsto ax^{\gamma}+b\mid a\in D,\, b\in\F,\, \gamma\in{\rm Gal}\,\F,\, \phi^{\gamma}=\phi\},\]
where now $\phi$ is an epimorphism $\F^*\to E$ for some group $E$ (necessarily cyclic) of order $d$. In addition, he considered those functions $f:\F\to\F$ for which there is a permutation $\kappa$ of $E$ such that
\begin{equation}\label{Lenstraeqn}
\phi(f(u)-f(v))=\kappa\phi(u-v) \; \hbox{for all} \; u\ne v\in \F.
\end{equation}
He showed that these functions $f$ form the normaliser $N(G(d))$ of $G(d)$ in the symmetric group on $\F$. In order to describe the elements of this group, let $K$ denote the subfield of $\F$ generated by $D$, and define a $K$-{\sl semilinear automorphism\/} of $\F$ to be an automorphism $\beta$ of the additive group of $\F$ for which there is a field automorphism $\gamma$ of $K$ satisfying $\beta(xy)=(\gamma x)(\beta y)$ for all $x\in K$ and $y\in\F$. Then $N(G(d))$ consists of the transformations $x\mapsto x^{\beta}+b$ of $\F$ such that $b\in\F$ and $\beta$ is a $K$-semilinear automorphism of $\F$.

As we have remarked, Carlitz's paper contains no references to graphs, groups or automorphisms, or to Paley. Some of these later generalisations by McConnel, Bruen, Levinger and Lenstra use group theory, to a varying extent, but none of them mentions graphs or Paley. The first proof to do that is the subject of the next section.


\subsection{Muzychuk's proof of the main theorem}

In 1987 Muzychuk~\cite{Muz87}  independently gave a full proof of Theorem~\ref{Paleyautos}. At that time he was a PhD student in Kiev, supervised by V.A.Ustimenko. This was towards the end of a long period during which contacts between Soviet mathematicians and those in the West were almost non-existent, so it is not surprising that the results of Carlitz and his successors were not known there, and were not cited in this paper; indeed, the only citations were to the Russian translations of the book by Cameron and van Lint~\cite{CvL80} for the definition (and name) of the Paley graphs, and of that by Serre~\cite{Ser} for some basic properties of finite fields. As in the case of some of the earlier proofs, the main argument involves the ingenious use of polynomials over finite fields. The paper was written in Russian, and published in a journal difficult to access outside the former Soviet Union, but as international contacts became much easier in the 1990s it became more widely known, with several recent citations listed in MathSciNet. An English translation, including an extension of the main theorem to cover cyclotomic schemes, is in preparation~\cite{Muz16}.

                                                                                                                                                                                                                                                                                                                                                                                                                                                                                                                                 
\subsection{Generalised Paley graphs} 

In 2009 Lim and Praeger~\cite{LP} introduced a class of graphs which generalise the Paley graphs, and in certain cases they found their automorphism groups. For consistency with earlier sections, we have changed their notation slightly.

Let $\F=\F_q$ for any prime power $q$, let $m$ be any divisor of $q-1$, and let $D$ be the unique subgroup of  order $m$ in $\F^*$. Like $P(q)$, a {\sl generalised Paley graph\/} $P=P^{(m)}(q)$ has vertex-set $\F=\F_q$, but with vertices $u$ and $v$ adjacent if and only if $u-v\in D$; in other words, $P$ is the Cayley graph for the additive group of $\F$, with connection set $D$. In order that $D=-D$, giving an undirected graph, we need to impose the following restriction:
\begin{itemize}
\item if $q$ is odd then $m$ is even.
\end{itemize}
Here, unlike Lim and Praeger, we will also assume that
\begin{itemize}
\item $D$ generates the additive group of $\F$,
\end{itemize}
so that $P$ is connected. For example, if $q\equiv 1$ mod~$(4)$ and $m=(q-1)/2$, then $P$ is the Paley graph $P(q)$.

It is clear that the transformations
\begin{equation}\label{genPaleyautos}
x\mapsto ax^{\gamma}+b\quad (a\in D,\, b\in\F,\, \gamma\in{\rm Gal}\,\F)
\end{equation}
of $\F$, which form a subgroup of index $d=(q-1)/m$ in $A\Gamma L_1(q)$, are all automorphisms of $P$. Unfortunately, the various extensions of Carlitz's theorem which we have discussed do not provide a converse, since if $m<(q-1)/2$ then condition (2) of Theorem~\ref{McCthm} is too restrictive: we need the weaker condition that if $\phi(u-v)=1$ then $\phi(f(u)-f(v))=1$. Similarly, Lenstra's condition in equation~(\ref{Lenstraeqn}) does not help in this situation.

However, Lim and Praeger~\cite{LP} proved a partial converse, showing that the transformations in~(\ref{genPaleyautos}) are the only automorphisms of $P$ provided $D$ is  `large' in the following sense:
\begin{itemize}
\item the index $d=|\F^*:D|$ divides $p-1$, where $q$ is a power of the prime $p$.
\end{itemize}
(This implies that $|D\cup\{0\}|>q/p$, so that $P$ is connected. However, there are examples where $P$ is connected, but $d$ does not divide $p-1$ and ${\rm Aut}\,P$ is larger: for instance, if $q=p^2$ and $m=2(p-1)$ then $P$ is a Hamming graph and ${\rm Aut}\,P$ is a wreath product $S_p\wr S_2$, of order $2(p!)^2$.)  Their proof of this partial converse (unlike the rest of the results in~\cite{LP}) depends on the classification of finite simple groups, and it would be interesting to find a more elementary proof.  For further properties of generalised Paley graphs, see~\cite{Jon13} for their regular surface embeddings, and~\cite{PP} for their product decompositions; see also~\cite[Remark~1.2]{PP} for applications of these graphs to topics such as Ramsey theory and synchronizing groups.


\subsection{The automorphism group of the Paley tournament}

The analogue of the Paley graph for a prime power $q\equiv 3$ mod~$(4)$ is a directed graph, called the {\sl Paley digraph}, or {\sl quadratic residue digraph}. This is a tournament, since every distinct pair of vertices are joined by a single arc.  The automorphism group of any finite tournament has odd order (since no element can transpose two vertices), so by the Feit-Thompson Theorem~\cite{FT} it is solvable. This makes the study of automorphism groups of tournaments relatively straightforward.

In 1970, Goldberg~\cite{Gol} used results on permutation groups to prove that the Paley digraph has automorphism group $A\Delta L_1(q)$. (There is no reference to Paley in this paper.) In a late note, the author wrote that his theorem was a special case of an unpublished result of Kantor, stated without proof in Dembowski's book~\cite[p.~98]{Dem} (see also Kantor's results on $2$-homogeneous groups~\cite{Kan72}, published in 1972). In 1972 Berggren~\cite{Ber} also proved this theorem, and showed that the Paley digraphs are the only  finite symmetric (vertex- and arc-transitive) tournaments.


\subsection{Automorphism groups of Hadamard matrices}\label{autoHadamard}

The automorphisms of a Hadamard matrix $H$ were defined by Hall~\cite{Hal62} to be the ordered pairs $(P,Q)$ of monomial matrices $P$ and $Q$, with non-zero entries $\pm 1$, such that $PHQ^T=H$; these form a group ${\rm Aut}\,H$. The element $(-I,-I)$ is a central involution $\sigma\in {\rm Aut}\,H$, and the quotient ${\overline{\rm Aut}\,H}={\rm Aut}\, H/\langle\sigma\rangle$ acts faithfully on the union of the sets of rows and columns of $H$.

The Paley-Hadamard matrices $H$ of type~I have order $m=q+1$ for prime powers $q\equiv 3$ mod~$(4)$. In~\cite{Hal62} Hall~showed that for these matrices, ${\overline{\rm Aut}\,H}$ contains $PSL_2(q)$ (not $P\Sigma L_2(q)$, as asserted in~\cite{Kan69b}, though it does indeed contain this group), acting on both the rows and the columns as a group of M\"obius transformations of the projective line $PG(1,q)={\mathbb P}^1(\F_q)=\F_q\cup\{\infty\}$. He also showed that if $q=11$ then ${\overline{\rm Aut}\,H}$ is strictly larger, acting as the Mathieu group $M_{12}$ on the rows and columns, whereas in 1969 Kantor~\cite{Kan69b} showed that  ${\overline{\rm Aut}\,H}=P\Sigma L_2(q)$ whenever $q>11$. 
By contrast, the automorphism group of the corresponding Paley design is $A\Delta L_1(q)$, a subgroup of index $q+1$ in $P\Sigma L_2(q)$, if $q\ge 19$, whereas it is $PSL_2(q)=P\Sigma L_2(q)$ if $q=7$ or $11$ (see Section~\ref{Todd}); the corresponding Paley graph $P(q)$ has automorphism group $A\Delta L_1(q)$ for all $q$ (see Section~\ref{charautos}).

The Paley-Hadamard matrices $H$ of type~II have order $m=2(q+1)$ for prime powers $q\equiv 1$ mod~$(4)$. If $q=5$ then $H$ is equivalent to the Paley-Hadamard matrix of type I and order $12$, discussed above. In 2008 De Launey and Stafford~\cite{deLS08} showed that if $q>5$ then ${\rm Aut}\,H$ has a subgroup of index $2$ isomorphic to $\Gamma L_2(q)/S$, where we identify $S$ with the group of scalar matrices $\lambda I\;(\lambda\in S)$ in $GL_2(q)$; the full group is obtained by adjoining an element of order $4$, whose square is the central involution in $\Gamma L_2(q)/S$, corresponding the matrices $\lambda I\;(\lambda\in\F^*\setminus S)$. Their proof uses the classification of finite simple groups, via the classification of 2-transitive finite permutation groups.



\section{Attribution and terminology}\label{Attrib}

We have seen that the papers by Paley~\cite{Pal} and Carlitz~\cite{Car}, frequently cited in connection with the Paley graphs and their automorphism groups, do not in fact mention graphs, groups or automorphisms. Indeed, inspection of their publication records suggests that neither of these mathematicians had much interest in either graph theory or group theory. This therefore raises the question of whether it is appropriate to refer to the graphs $P(q)$ as `Paley graphs', or to attribute Theorem~\ref{Paleyautos}, describing their automorphism groups, to Carlitz.

In view of our remarks in the preceding sections, to refer to Theorem~\ref{Paleyautos} as `Carlitz's Theorem' seems slightly over-generous, even though he did the hard work in providing most of the argument for the difficult half of the proof. (In any case, Carlitz is hardly short of recognition: MathSciNet lists 732 publications by him, together with (at the time of writing) 
2308 citations of his work, and 260 publications with his name in the title.) Perhaps it would be more correct to reserve this term for the result he actually proved in~\cite{Car}, namely Theorem~\ref{Carlitzthm}, and to refer to Theorem~\ref{Paleyautos} as a straightforward corollary to his theorem. 

 At first sight, Paley's connection with the graphs $P(q)$ seems to be rather tenuous. As was later shown, Hadamard matrices, including those constructed by him, give rise to graphs, in the sense that normalising the matrix, deleting the first row and column, and replacing each entry $-1$ with $0$, produces the adjacency matrix of a graph. Paley's first construction yields directed graphs of prime power order $q\equiv 3$ mod~$(4)$ (the Paley tournaments), rather than the undirected graphs of order $q\equiv 1$ mod~$(4)$ which bear his name, while his second construction yields undirected graphs of order $2q+1$, where $q\equiv 1$ mod~$(4)$. However, one ingredient of this second construction~\cite[p.~314]{Pal} is a (non-Hadamard) matrix $(B_{ij})$ of order $q+1$ such that deleting its first row and column yields the Jacobsthal matrix $Q=(\chi(j-i))$ for $\F_q$, and hence (after replacing each entry $-1$ with $0$) the adjacency matrix for the graph $P(q)$. In this sense, the Paley graphs do arise naturally from Paley's paper, though not directly from the Hadamard matrices he constructed.
 
Whether this justifies naming these graphs after Paley is debatable. The first person to describe these graphs in the literature seems to have been Sachs~\cite{Sac}, in 1962, who restricted attention to prime values of $q$, followed independently in 1963 by Erd\H os and R\'enyi~\cite{ER}, who considered the general case. In 1971 Higman~\cite{Hig} constructed the Paley graphs as examples of strongly regular graphs, but did not name them or refer to Paley; he referred to certain rank 3 permutation groups acting on $P(q)$ as `of Singer type', a terminology which does not seem to have survived. 

The main link between Paley and the graphs $P(q)$ is the use of quadratic residues, and in particular their combinatorial property given by equation~(\ref{chisum}). However, this result can be traced back at least to Jacobsthal~\cite{Jac06, Jac07}, a generation earlier. Whoever coined the term `Paley graph' probably did so as shorthand for a more accurate but clumsy phrase such as `graph based on Paley's construction', rather than as a deliberate attribution. As suggested to the author by Mikhail Muzychuk, the terms `Paley graph' and `Carlitz's theorem' appear to be further instances of Stigler's Law of Eponymy (which is, of course, itself also a misattribution~\cite{Stig}).

\section{Acknowledgements}

The author is grateful to Brian Alspach, Norman Biggs, B\'ela Bollob\'as, Andries Brouwer, Peter Cameron, Chris Godsil, Willem Haemers, Joshua Insley, Mikhail  Muzychuk, Cheryl Praeger and Don Taylor for helpful comments, information  and suggestions, and in particular to Mikhail Klin and Franz Lemmermeyer for their valuable advice on the extensive literature related to the Paley graphs and on the mathematical works of Gauss.

.

\bigskip

\noindent School of Mathematics

\noindent University of Southampton

\noindent Southampton SO17 1BJ

\noindent UK

\medskip

\noindent G.A.Jones@maths.soton.ac.uk

\end{document}